\documentclass[reqno,12pt,a4paper]{amsart}

\usepackage{mathrsfs}
\usepackage{amssymb}
\usepackage{amsfonts}
\usepackage{latexsym}
\usepackage{amsthm}

\usepackage{graphicx}
\def\lb{\label}

\newcommand{\er}[1]{\textrm{(\ref{#1})}}

\begin{document}

%%%%%%%%%% Some definitions %%%%%%%%%%

%%%%%%%% Equations, theorems %%%%%%%%%
\renewcommand{\theequation}{\arabic{section}.\arabic{equation}}
\theoremstyle{plain}
\newtheorem{theorem}{\bf Theorem}[section]
\newtheorem{lemma}[theorem]{\bf Lemma}
\newtheorem{corollary}[theorem]{\bf Corollary}
\newtheorem{proposition}[theorem]{\bf Proposition}
\newtheorem{definition}[theorem]{\bf Definition}
\newtheorem{remark}[theorem]{\it Remark}
%\theoremstyle{remark}
%\newtheorem{remark}[theorem]{\bf Remark}

%%%%% Alphabet %%%%%
\def\a{\alpha}  \def\cA{{\mathcal A}}     \def\bA{{\bf A}}  \def\mA{{\mathscr A}}
\def\b{\beta}   \def\cB{{\mathcal B}}     \def\bB{{\bf B}}  \def\mB{{\mathscr B}}
\def\g{\gamma}  \def\cC{{\mathcal C}}     \def\bC{{\bf C}}  \def\mC{{\mathscr C}}
\def\G{\Gamma}  \def\cD{{\mathcal D}}     \def\bD{{\bf D}}  \def\mD{{\mathscr D}}
\def\d{\delta}  \def\cE{{\mathcal E}}     \def\bE{{\bf E}}  \def\mE{{\mathscr E}}
\def\D{\Delta}  \def\cF{{\mathcal F}}     \def\bF{{\bf F}}  \def\mF{{\mathscr F}}
\def\c{\chi}    \def\cG{{\mathcal G}}     \def\bG{{\bf G}}  \def\mG{{\mathscr G}}
\def\z{\zeta}   \def\cH{{\mathcal H}}     \def\bH{{\bf H}}  \def\mH{{\mathscr H}}
\def\e{\eta}    \def\cI{{\mathcal I}}     \def\bI{{\bf I}}  \def\mI{{\mathscr I}}
\def\p{\psi}    \def\cJ{{\mathcal J}}     \def\bJ{{\bf J}}  \def\mJ{{\mathscr J}}
\def\vT{\Theta} \def\cK{{\mathcal K}}     \def\bK{{\bf K}}  \def\mK{{\mathscr K}}
\def\k{\kappa}  \def\cL{{\mathcal L}}     \def\bL{{\bf L}}  \def\mL{{\mathscr L}}
\def\l{\lambda} \def\cM{{\mathcal M}}     \def\bM{{\bf M}}  \def\mM{{\mathscr M}}
\def\L{\Lambda} \def\cN{{\mathcal N}}     \def\bN{{\bf N}}  \def\mN{{\mathscr N}}
\def\m{\mu}     \def\cO{{\mathcal O}}     \def\bO{{\bf O}}  \def\mO{{\mathscr O}}
\def\n{\nu}     \def\cP{{\mathcal P}}     \def\bP{{\bf P}}  \def\mP{{\mathscr P}}
\def\r{\rho}    \def\cQ{{\mathcal Q}}     \def\bQ{{\bf Q}}  \def\mQ{{\mathscr Q}}
\def\s{\sigma}  \def\cR{{\mathcal R}}     \def\bR{{\bf R}}  \def\mR{{\mathscr R}}
                \def\cS{{\mathcal S}}     \def\bS{{\bf S}}  \def\mS{{\mathscr S}}
\def\t{\tau}    \def\cT{{\mathcal T}}     \def\bT{{\bf T}}  \def\mT{{\mathscr T}}
\def\f{\phi}    \def\cU{{\mathcal U}}     \def\bU{{\bf U}}  \def\mU{{\mathscr U}}
\def\F{\Phi}    \def\cV{{\mathcal V}}     \def\bV{{\bf V}}  \def\mV{{\mathscr V}}
\def\P{\Psi}    \def\cW{{\mathcal W}}     \def\bW{{\bf W}}  \def\mW{{\mathscr W}}
\def\o{\omega}  \def\cX{{\mathcal X}}     \def\bX{{\bf X}}  \def\mX{{\mathscr X}}
\def\x{\xi}     \def\cY{{\mathcal Y}}     \def\bY{{\bf Y}}  \def\mY{{\mathscr Y}}
\def\X{\Xi}     \def\cZ{{\mathcal Z}}     \def\bZ{{\bf Z}}  \def\mZ{{\mathscr Z}}
\def\O{\Omega}

\def\mb{{\mathscr b}}
\def\mh{{\mathscr h}}
\def\me{{\mathscr e}}
\def\mk{{\mathscr k}}
\def\mz{{\mathscr z}}
\def\mx{{\mathscr x}}

\newcommand{\gA}{\mathfrak{A}}          \newcommand{\ga}{\mathfrak{a}}
\newcommand{\gB}{\mathfrak{B}}          \newcommand{\gb}{\mathfrak{b}}
\newcommand{\gC}{\mathfrak{C}}          \newcommand{\gc}{\mathfrak{c}}
\newcommand{\gD}{\mathfrak{D}}          \newcommand{\gd}{\mathfrak{d}}
\newcommand{\gE}{\mathfrak{E}}
\newcommand{\gF}{\mathfrak{F}}           \newcommand{\gf}{\mathfrak{f}}
\newcommand{\gG}{\mathfrak{G}}           %\newcommand{\gg}{\mathfrak{g}}
\newcommand{\gH}{\mathfrak{H}}           \newcommand{\gh}{\mathfrak{h}}
\newcommand{\gI}{\mathfrak{I}}           \newcommand{\gi}{\mathfrak{i}}
\newcommand{\gJ}{\mathfrak{J}}           \newcommand{\gj}{\mathfrak{j}}
\newcommand{\gK}{\mathfrak{K}}            \newcommand{\gk}{\mathfrak{k}}
\newcommand{\gL}{\mathfrak{L}}            \newcommand{\gl}{\mathfrak{l}}
\newcommand{\gM}{\mathfrak{M}}            \newcommand{\gm}{\mathfrak{m}}
\newcommand{\gN}{\mathfrak{N}}            \newcommand{\gn}{\mathfrak{n}}
\newcommand{\gO}{\mathfrak{O}}
\newcommand{\gP}{\mathfrak{P}}             \newcommand{\gp}{\mathfrak{p}}
\newcommand{\gQ}{\mathfrak{Q}}             \newcommand{\gq}{\mathfrak{q}}
\newcommand{\gR}{\mathfrak{R}}             \newcommand{\gr}{\mathfrak{r}}
\newcommand{\gS}{\mathfrak{S}}              \newcommand{\gs}{\mathfrak{s}}
\newcommand{\gT}{\mathfrak{T}}             \newcommand{\gt}{\mathfrak{t}}
\newcommand{\gU}{\mathfrak{U}}             \newcommand{\gu}{\mathfrak{u}}
\newcommand{\gV}{\mathfrak{V}}             \newcommand{\gv}{\mathfrak{v}}
\newcommand{\gW}{\mathfrak{W}}             \newcommand{\gw}{\mathfrak{w}}
\newcommand{\gX}{\mathfrak{X}}               \newcommand{\gx}{\mathfrak{x}}
\newcommand{\gY}{\mathfrak{Y}}              \newcommand{\gy}{\mathfrak{y}}
\newcommand{\gZ}{\mathfrak{Z}}             \newcommand{\gz}{\mathfrak{z}}

\def\ve{\varepsilon} \def\vt{\vartheta} \def\vp{\varphi}  \def\vk{\varkappa} \def\vr{\varrho}

\newcommand{\dA}{\mathbb A}
\def\B{{\mathbb B}}
\def\C{{\mathbb C}}
\def\dD{{\mathbb D}}
 \def\dE{\mathbb E}
  \def\dG{\mathbb G}
   \def\dF{\mathbb F}
 \def\dJ{{\mathbb J}}
 \def\dI{{\mathbb I}}
\def\K{{\mathbb K}}
\def\N{{\mathbb N}}
 \def\R{{\mathbb R}}
 \def\S{{\mathbb S}}
 \def\T{{\mathbb T}}
 \def\Z{{\mathbb Z}}

%%%%% Arrows %%%%%

\def\la{\leftarrow}              \def\ra{\rightarrow}     \def\Ra{\Rightarrow}
\def\ua{\uparrow}                \def\da{\downarrow}
\def\lra{\leftrightarrow}        \def\Lra{\Leftrightarrow}
\newcommand{\abs}[1]{\lvert#1\rvert}
\newcommand{\br}[1]{\left(#1\right)}

\def\lan{\langle} \def\ran{\rangle}

%%%%% Typography %%%%%

\def\lt{\biggl}                  \def\rt{\biggr}
\def\ol{\overline}               \def\wt{\widetilde}
\def\no{\noindent}

%%%%% Math signs %%%%%

\let\ge\geqslant                 \let\le\leqslant
\def\lan{\langle}                \def\ran{\rangle}
\def\/{\over}                    \def\iy{\infty}
\def\sm{\setminus}               \def\es{\emptyset}
\def\ss{\subset}                 \def\ts{\times}
\def\pa{\partial}                \def\os{\oplus}
\def\om{\ominus}                 \def\ev{\equiv}
\def\iint{\int\!\!\!\int}        \def\iintt{\mathop{\int\!\!\int\!\!\dots\!\!\int}\limits}
\def\el2{\ell^{\,2}}             \def\1{1\!\!1}
\def\sh{\sharp}
\def\wh{\widehat}
\def\bs{\backslash}
\def\na{\nabla}
%%%%% Math operations %%%%%

\def\sh{\mathop{\mathrm{sh}}\nolimits}
\def\all{\mathop{\mathrm{all}}\nolimits}
\def\Area{\mathop{\mathrm{Area}}\nolimits}
\def\arg{\mathop{\mathrm{arg}}\nolimits}
\def\const{\mathop{\mathrm{const}}\nolimits}
\def\det{\mathop{\mathrm{det}}\nolimits}
\def\diag{\mathop{\mathrm{diag}}\nolimits}
\def\diam{\mathop{\mathrm{diam}}\nolimits}
\def\dim{\mathop{\mathrm{dim}}\nolimits}
\def\dist{\mathop{\mathrm{dist}}\nolimits}
\def\Im{\mathop{\mathrm{Im}}\nolimits}
\def\Iso{\mathop{\mathrm{Iso}}\nolimits}
\def\Ker{\mathop{\mathrm{Ker}}\nolimits}
\def\Lip{\mathop{\mathrm{Lip}}\nolimits}
\def\rank{\mathop{\mathrm{rank}}\limits}
\def\Ran{\mathop{\mathrm{Ran}}\nolimits}
\def\Re{\mathop{\mathrm{Re}}\nolimits}
\def\Res{\mathop{\mathrm{Res}}\nolimits}
\def\res{\mathop{\mathrm{res}}\limits}
\def\sign{\mathop{\mathrm{sign}}\nolimits}
\def\span{\mathop{\mathrm{span}}\nolimits}
\def\supp{\mathop{\mathrm{supp}}\nolimits}
\def\Tr{\mathop{\mathrm{Tr}}\nolimits}
\def\BBox{\hspace{1mm}\vrule height6pt width5.5pt depth0pt \hspace{6pt}}
\def\where{\mathop{\mathrm{where}}\nolimits}
\def\as{\mathop{\mathrm{as}}\nolimits}

%%%%%%%%%%%%% specialities %%%%%%%%%%%%%%

\newcommand\nh[2]{\widehat{#1}\vphantom{#1}^{(#2)}}
%{{\mathop{#1}\limits^\wedge}\vphantom{#1}^{(#2)}}
\def\dia{\diamond}

\def\Oplus{\bigoplus\nolimits}

%%%%%%%%%%% End of definitions %%%%%%%%%%

%%%%% OLD OLD OLD

\def\qqq{\qquad}
\def\qq{\quad}
\let\ge\geqslant
\let\le\leqslant
\let\geq\geqslant
\let\leq\leqslant
\newcommand{\ca}{\begin{cases}}
\newcommand{\ac}{\end{cases}}
\newcommand{\ma}{\begin{pmatrix}}
\newcommand{\am}{\end{pmatrix}}
\renewcommand{\[}{\begin{equation}}
\renewcommand{\]}{\end{equation}}
\def\eq{\begin{equation}}
\def\qe{\end{equation}}
\def\[{\begin{equation}}
\def\bu{\bullet}

\title[{Eigenvalues of Schr\"odinger operators  on finite and infinite
intervals }] {Eigenvalues of Schr\"odinger operators  on finite
and infinite intervals}

\date{\today}
\author[Evgeny Korotyaev]{Evgeny L. Korotyaev}
\address{Department of Math. Analysis, Saint-Petersburg State University,
Universitetskaya nab. 7/9, St.
Petersburg, 199034, Russia, \ korotyaev@gmail.com, \
e.korotyaev@spbu.ru}


\begin{abstract}
We consider a Sturm-Liouville operator a with integrable potential $q$ on the unit interval $I=[0,1]$.
We consider a Schr\"odinger operator with a real compactly supported
potential on the half line and on the line, where this  potential coincides with $q$ on the unit  interval  and vanishes outside $I$.
We determine the relationships between eigenvalues of such operators and obtain estimates of eigenvalues in terms of potentials.
\end{abstract}

\maketitle

%\noindent{\small \textbf{Keywords:} eigenvalues, estimates} \ \ \ \
%{\bf Preliminary version}a

%\noindent {\small \textbf{AMS Subject classification:} }

\vskip 0.25cm
\section {Introduction and main results}
\setcounter{equation}{0}

%\subsection{Introduction}

In this paper we discuss eigenvalues of  Schr\"odinger operators defined by
\[
\lb{defT}
\begin{aligned}
&Case \ 1: \qqq\textstyle  T y=-y''+qy \qq {\rm on}\qq L^2(\R_+),\qq
{\rm with} \ y(0)=0, \\
&Case \ 2: \qqq \textstyle \wt T y=-y''+qy\qq {\rm on}\qq L^2(\R_+),\qq {\rm with}\
\ y'(0)=0,\\
&Case \ 3: \qqq \textstyle \cT y=-y''+qy\qq {\rm on}\qq L^2(\R ),
\end{aligned}
\]
and defined by \er{H01}-\er{per} on the unit interval.
We assume that the potential $q$ satisfies
\[
\lb{dq}
\supp q\in [0,1],\qqq q\in L_{real}^1(0,1).
\]
 It is well known (\cite{F59}, \cite{DT79}, \cite{M86}) that their spectrum
consists of an absolutely continuous part $[0,\iy)$ plus a finite
number of simple negative eigenvalues given by
$$
\begin{aligned}
\s_{ac}(T)=[0,\iy),\qq
\s_{d}(T)=\{E_1<\dots <E_m<0\},
\\
\s_{ac}(\wt T)=[0,\iy),\qq
\s_{d}(\wt T)=\{\wt E_1<\dots <\wt E_N<0\},
\\
\s_{ac}(\cT)=[0,\iy),\qq  \s_{d}(\cT)=\{\cE_1<\dots <\cE_\cN<0\}.
\end{aligned}
$$
We  introduce Sturm-Liouville operators on the interval $[0,1]$ with the Dirichlet and
Neumann boundary conditions:
\[
\lb{H01}
\begin{aligned}
H_0f=-f''+qf,\qqq f(0)=f(1)=0,\\
H_{1} f=-f''+qf,\qqq f'(0)=f'(1)=0,
\end{aligned}
\]
and with so-called mixed boundary conditions:
\[
\lb{H0110}
\begin{aligned}
H_{01} f=-f''+qf, \qqq f(0)=f'(1)=0,
\\
H_{10} f=-f''+qf, \qqq f'(0)=f(1)=0,
\end{aligned}
\]
and the operator $H_\pi$ on $[0,2]$ with 2-periodic boundary conditions:
\[
\lb{per}
\begin{aligned}
H_\pi f=-f''+q_\pi f, \qqq q_\pi=\ca q(x), \  & x\in [0,1]\\
                                     q(x-1), \ & x\in [1,2]\ac .
\end{aligned}
\]

Let
$\m_n$ and $\n_0,  \n_n, n\ge 1$ be eigenvalues of $H_0$ and $H_{1}$ respectively.
Let $\t_n$ and $\vr_n, n\ge 1$ be eigenvalues of $ H_{01}$ and
$H_{10}$ respectively. All these eigenvalues are simple.
It is well known, that the spectrum of $H_\pi$ discrete and its  eigenvalues
$\l_0^+, \l_n^\pm, n\ge 1$ satisfy $\l^+_{n-1}<\l^-_n \leq  \l^+_{n},\ n\geq 1$.
Moreover, all these eigenvalues satisfy
\[
\lb{bx}
\begin{aligned}
& \n_0\le \l_0^+<  \ol {\vr_1, \t_1}<\l_1^- \le \ol {\m_1,\n_1}\le \l_1^+<\ol {\vr_2, \t_2}< \l_2^- \le \ol {\m_2,\n_2}\le \l_2^+<...,\\
& \l_n^\pm , \n_n,\m_n=(\pi n)^2+c_0+o(1),\qq \vr_n,\t_n=\pi^2     \textstyle  (n-{1\/2})^2+c_0+o(1)\qq \as \ n\to \iy,
\end{aligned}
\]
where $\ol {u,v}$ denotes $\min \{u,v\}\le \max \{u,v\}$ for shortness and $c_0=\int_0^1qdx$.

Given a self-adjoint, bounded from below operator $A$ whose negative
spectrum is discrete, we denote by $n_-(A)$ the number of its
negative eigenvalues, counted with multiplicities.

\begin{theorem}
\lb{T1} i) The eigenvalues $E_1<\dots <E_m<0$  of $T$ for
the Case 1 satisfy
\[
\lb{D1} n_-(T)=n_-(H_{01})\ge n_-(H_{0}),
\]
%and
\[
\lb{D2}
\begin{aligned}
\t_1<E_1<\m_1<\t_2<E_2<\m_2<\dots <\t_m<E_m<\min \{0,\m_m\}.
%\m_j<r_{j,1}^2\le...\le r_{j,s_j}^2<\t_{j+1}.
\end{aligned}
\]
ii) The eigenvalues $\wt E_1<\dots <\wt E_N<0$ of $\wt T$ for the Case 2
satisfy
\[
\lb{N1}  n_-(\wt T)=n_-(H_1),
\]
and
\[
\lb{N2} \n_0<\wt E_1<\vr_1<\n_1<\wt E_2<\vr_2<\dots
<\n_{N-1}<\wt E_N<\min \{0,\vr_N\},
\]
\[
\lb{N3}
\begin{aligned}
 \wt E_1< E_1<\wt E_2<E_2<\dots <\wt E_N,\qqq \ca   {\rm  if} \ \ \t_N\ge 0 \Rightarrow \ m=N-1\\
                                                      {\rm  if} \ \ \t_N<0 \Rightarrow \ m=N \ \  {\rm and} \ \ \wt E_N<E_N\ac.
\end{aligned}
\]

\no iii) The operator $\cT$ for the Case 3 has the eigenvalues $\cE_1<....<\cE_N$, which  satisfy
\[
\lb{N1x} n_-(\cT)=n_-(\wt T)=n_-(H_1),
\]
\[
\lb{N2x} \n_0<\wt E_1<\cE_1<E_1< \n_1 <\wt E_2<\cE_2<\dots <
\n_{N-1}<\wt E_{N}<\cE_N<\min \{0,E_N\}.
\]
\end{theorem}

{\bf Remark.} Thus eigenvalues $\{\t_n\}$ control $\{E_n\}$; $\{\n_n\}$
 control both $\{\wt E_n\}$ and $\{\cE_n\}$.

We describe eigenvalues for $\cT_e y=-y''+q_e y$ on $L^2(\R)$ with even potentials $q_e(x)=q(|x|)$ for all $x\in \R$.
Using Theorem \ref{T1} and the known identity \er{tw} we obtain
%Define the space $L_{even}^1(0,1)=\{q\in L^1(0,1): q=q(1-x)\  \forall \ x\in (0,1)\}$.

\begin{corollary}
\lb{T2} The eigenvalues   $\cE_1^e<\dots <\cE_{n_e}^e<0$  of $\cT_e$
coincide with eigenvalues of the operators $T$ and $\wt T$ given by
\er{defT} and satisfy
\[
\lb{ev1}
\s_{d}(\cT_e)=\s_{d}(T)\cup \s_{d}(\wt T),\qq
n_e=n_-(\cT_e)=N+m=\ca  2N-1 \ \  {\rm  if} \ \ \t_N\ge 0\\
          2N  \ \ \  \ \ \ \  {\rm  if} \ \ \t_N<0\ac,
\]
\[
\lb{ev2}
\begin{aligned}
\n_0<\cE_1^e <\vr_1; \qq
\t_1<\cE_2^e<\m_1;\qq \n_1<\cE_3^e <\vr_2; \qq \t_2<\cE_4^e<\m_2;...
\end{aligned}
\]
\end{corollary}

%\newpage

{\bf Localization of resonances.} We discuss  resonances for the
class compactly supported potentials $\cP=\big\{f\in
L_{real}^1(\R_+): \sup\supp f=1\big\}$. We define the Jost solutions
$f_{+}(x,k)$ of the equation
\[
\lb{e1a}
  -f''+q(x)f=k^2f,\ \ \  \ \ \ k\in \C \sm \{0\},
\]
under the conditions: $f_+(x,k )=e^{ixk }$ for $x>1$.  The Jost
function $\p=f_+(0,\cdot)$ is entire and satisfies
\[
\lb{asj} \p(k)=1+ O(1/k)\qqq as \qq  |k|\to \iy,\qq k\in \ol\C_+,
\]
uniformly in $\arg k\in [0,\pi]$.
%The function $\p(k)$ is real on the imaginary line $\Re k =0$, and
%then $\ol \p(k )=\p(-k )$ for all $ k \in \R$.
The function $\p(k)$ has $m$ simple zeros $k_j, j\in
\N_m=\{1,..,m\}$ in $\C_+$, given by $ k_{j}=i|E_{j}|^{1\/2}, j\in
\N_m$,  possibly one simple zero at 0, and an infinite number
(so-called resonances) in $\C_-$: $ 0\le |k_{m+1}|\le |k_{m+2}|\le
....... $. The function $\p$ has an odd number $\geq 1$ of zeros on
any interval $I_j:=(-k_{j}, -k_{j+1}), j\in \N_{m-1}$ and an even
number $\geq 0$ of zeros on the interval $I_m:=(-i|k_{m}|,0]$
counted with multiplicity.
%Define the intervals $I_j=(-k_{j},-k_{j+1}), j\in \N_{m-1},$ and
%$(-k_{m},0]$ on  $i\R_-$. The zeros of function $\p$ on $I_j$ satisfies: each number
%$\#(I_j)\ge 1, j\in \N_{m-1}$ is odd and the number $\#(I_m)\ge 0$
%is  even.
%\[ \#(I_j)=\gn_j, \qq \ca \gn_j\ge 1 \ {\rm is \ odd \ if} \ j\in \N_{m-1}\\
%                     \gn_m\ge 0 \ {\rm }\ac \]
Let  $\#(E)$ be  the number of zeros of $\p$ (counted with
multiplicity) in the set $E\ss \C$.
 We describe the localization of resonances on the
interval $[-k_1,0]$ in terms of eigenvalues .

\begin{theorem}
\lb{T1x} i) Let $q\in \cP$ and let a resonance   $k_o\in I_j\ss
i\R_-$ for some $j\in \N_{m}$. Then
\[
\lb{D2r}
\begin{aligned}
\m_j<k_o^2<\t_{j+1},\qq j\in \N_{m}=\{1,..,m\}.
\end{aligned}
\]
ii) For any finite sequence $p_1, p_2,..,p_m\in i\R_+$ (labeled by
$|p_1|>|p_2|>...>|p_m|>0$)  and any odd integers $s_1,...,s_{m-1}\in
\N$ and even integer $s_m\ge 0$ there exists a potential $q\in \cP$,
such that each its eigenvalue $E_j=p_j^2<0$ and $\# (I_j)=s_j,
j=\N_m$.

\end{theorem}

{\bf Remark.} 1) Thus eigenvalues $\m_j, \t_{j+1}$ control
resonances on the interval  $I_j$.

2) If $\m_m\ge 0$, then there are no resonances on the interval
$(-k_m, 0]$.

3) Jost and Kohn \cite{JK53}  adapted the method of Gelfand and
Levitan to determine a potential (exponentially decaying) of a
Schr\"odinger operator by the spectral data. Theorem \ref{T1x} is
used in the paper \cite{K19} to show that  Jost and Kohn solution is
not complete and there exists another solution, which gives a
compactly supported potential.

%\newpage

{\bf Estimates.} Recall  the Bargmann inequality \cite{B52}: let  in
general, $xq(x)$ belong to $L^1(0,\iy)$, then
\[
\lb{Be}
n_-(T)\le \int_0^\iy xq_-(x)dx,
\]
where $q_-=\max \{-q, 0\}$. The  case of distributions was discussed
 in \cite{AKM10}, and more operators was studied in \cite{AKMN13}.  Recall the Calogero-Cohn inequality \cite{Ca65},
\cite{Co65}: if $q, |q|^{1\/2}\in L^1(0,\iy)$ and if $q\le 0$ is
monotone, then
\[
\lb{CC}
n_-(T)\le {2\/\pi}\int_0^\iy |q(x)|^{1\/2}dx.
\]

We recall the estimates from \cite{K98}, \cite{K00}: Consider the
operator $H_\pi$ with the potential $q\in \mH_0=\{q\in
L_{real}^2(0,1): \int_0^1 q(x)dx=0\} $ having eigenvalues $\l_0^+,
\l_n^\pm, n\ge 1$. Let $\g_n=\l_n^+-\l_n^-\ge 0, n\ge 1$ and
$\g=(\sum_{n\ge 1}\g_n^2)^{1\/2}\ge 0$. Then
 \[
\lb{esg1}
\begin{aligned}
\|q\|\le 2\g\max \{1, \g^{1\/3}\},
\\
\g\leq 2\|q\|\max \{1, \|q\|^{1\/3}\}.
\end{aligned}
\]

\begin{corollary}
\lb{T3} Let the operators $ H_{0}, H_{01}, H_{\pi}$ given by
\er{H01}-\er{per}.

i) Let $\gS=\sum_{\m_n<0}|\m_n|^{1\/2}$ or
$\gS={1\/2}\sum_{n>0,\l_n^\pm<0}|\l_n^\pm|^{1\/2}$.
 Then
 \[
\lb{nu1} {\textstyle {1\/2}}(n_-(H_{\pi})-1)\le n_-(H_{0})\le  n_-(H_{01})\le \int_0^1 xq_-(x)dx,
\]
 \[
\lb{nu2} \gS\le {1\/2}\int_0^1 q_-(x)dx,
\]
 Let $q(x)$ be monotone. Then
 \[
\lb{nu3} n_-(H_{0})\le n_-(H_{01})\le {2\/\pi}\int_0^1
q_-(x)^{{1\/2}}dx.
\]
ii) Let $q\in L^2(0,1)$ and $\int_0^1 q(x)dx=0$. Then
 \[
\lb{nu2x}
\begin{aligned}
& {\textstyle {1\/2}}(n_-(H_{\pi})-1)\le n_-(H_{0})\le
n_-(H_{01})\le 2\g\max \{1, \g^{1\/3}\},
\\
& \gS\le 2\g\max \{1, \g^{1\/3}\}.
\end{aligned}
\]

\end{corollary}

{\bf Remark.} Instead of $q$ in Corollary \ref{T3} we can use $q-\ve
$ for $\ve\in \R$

There are a lot of results devoted to estimates of negative eigenvalues of
$\dim d=1$ Schr\"odinger operators, see \cite{Ca65}, \cite{Co65}, \cite{B52}, \cite{DHS03}, \cite{HLT98},
\cite{S13},   \cite{W96} and references therein.

There exist many results about Sturm-Liouville operators on the unit
interval, see  \cite{DT84}, \cite{IT83}, \cite{KC09}, \cite{M86},
\cite{PT87} and on the circle, see \cite{GT87}, \cite{K99},
\cite{M86} and references therein. Unfortunately there only few
results about estimates for such Sturm-Liouville operators, for
example,  two-sided estimates of periodic potentials in terms of
gap-lengths see \er{esg1}.

\section {Proof}
\setcounter{equation}{0}

\subsection {Fundamental solutions}
Let $\vp(x,\l),  \vt(x,\l)$ be the solutions of the equation
\[
\lb{Fu1} -y''+qy=\l y, \ \ \ \l\in \C ,
\]
under the conditions
$\vp'(0,\l)=\vt(0,\l)=1$  and $\vp(0,\l)=\vt'(0,\l)=0$.
If $q=0$, then the corresponding  fundamental solutions are given by $\vp_0(x,\l)={\sin \sqrt \l x\/\sqrt \l}$ and $ \vt_0(x,\l)=\cos x\sqrt \l $. Recall the well-known results.

%\bigskip

\begin{lemma}
\lb{TFu21}  Let $q\in L^1(0,1)$. Then
the functions $\vp(1,\l),\vp'(1,\l)$, $\vt(1,\l),
\vt'(1,\l)$ are entire and satisfy:
\[
\begin{aligned}
\lb{Fu27} \textstyle \vp(1,\l)-{\sin \sqrt \l \/\sqrt \l}=e^{|\Im \sqrt \l|}O(\l^{-1}),
\\
\vp'(1,\l)-\cos \sqrt \l = e^{|\Im \sqrt \l|}O(\l^{-{1\/2}}),
\end{aligned}
\]
and
\[
\begin{aligned}
\lb{Fu27x} \vt(1,\l)-\cos \sqrt \l =e^{|\Im \sqrt \l|}O(\l^{-{1\/2}}),
\\
\vt'(1,z)+\sqrt \l\sin \sqrt \l =e^{|\Im \sqrt \l|}O(1),
\end{aligned}
\]
as $|\l|\to \iy$, uniformly in $\arg \l\in [0,2\pi]$, and, in
particular,
\[
\lb{aff1}
\vp(1,\l), \vp'(1,\l), \vt(1,\l), \vt'(1,\l)\to +\iy\qq  as \qq \l\to -\iy.
\]
\end{lemma}

Note that $\{\m_n\}, \{\n_n\}, \{\t_n\}, \{\vr_n\}$ are the zeros of $\vp
(1,\l), \vt'(1,\l), \vp'(1,\l), \vt(1,\l)$ respectively.

\subsection{The periodic case}

%For $q$ satisfying \er{dq} we define the 1-periodic  potential $q_\pi=q(x), x\in [0,1]$.
We consider the 2-periodic operator $H_\pi= -{d^2\/dx^2}+q_\pi(x)$  on
$[0,2]$ given by \er{per}. It is well known, that the spectrum of $H_\pi$ discrete and is eigenvalues
$\l_0^+, \l_n^\pm, n\ge 1$, which satisfy $\l^+_{n-1}<\l^-_n \leq  \l^+_{n},\ n\geq 1$.
% absolutely
%continuous, which are union  of bands $\S_n=[A^+_{n-1},A^-_n ],$
%where $A^+_{n-1}< A^-_n \leq  A^+_{n},\ n\geq 1$.
%The eigenvalues $\l_{n}^+$ and $\l_{n+1}^-$  are separated by the gap $\g_n=(\l^-_n,
%\l^+_n )$.
% If the gap degenerates, that is $\g_n=\es$,
%then the corresponding eigenvalues $\l_{n}^+$ and $\l_{n+1}^-$ merge.
The sequence $\l_0^+<\l_1^-\leq \ \l_1^+\ <\dots$ is the spectrum of
equation
\[
\lb{ee}
-y''+q_\pi y=\l y
\]
 with 2-periodic boundary conditions, i.e.
$y(x+2)=y(x), x\in \R$. If $\l_n^-=\l_n^+$ for some $n$, then this
number $\l_n^{\pm}$ is a double eigenvalue of Eq. \er{ee} with
2-periodic  boundary conditions. The lowest  eigenvalue $\l_0^+$ is
always simple, and the corresponding eigenfunction is 1-periodic.
The eigenfunctions corresponding to the eigenvalue $\l_n^{\pm}$ are
1-periodic, when $n$ is even and they are antiperiodic, i.e.
$y(x+1)=-y(x),\ \ x\in\R$, when $n$ is odd. We introduce the Lyapunov function
$\D={1\/2}(\vp'(1,\l)+\vt(1,\l))$
 It is
well known that
\[
\lb{3.9}
\D (\l_0^+)=1,\qq
\D (\l_n^{\pm})=(-1)^n,  \qqq (-1)^n\D ([\l_n^-,\l_n^+])\ss [1,\iy), \ n\geq 1.
\]
%Here and below we use the notation $(\ ') =\pa /\pa x ,\ \ \dot
%{}=\pa /\pa \l$.
\subsection{Jost functions}
We recall well-known result about the Jost function from
\cite{F59}, \cite{M86}, \cite{DT79}. The Schr\"odinger equation
\[
\lb{e1a}
  -f''+q(x)f=k^2f,\ \ \  \ \ \ k\in \C \sm \{0\},
\]
has unique solutions $f_{\pm}(x,k)$  such that $f_+(x,k )=e^{ixk }$
for $x>1$ and $f_-(x,k)=e^{-ik x}$ for  $x<0$. Outside the support
of $q$ any solutions of \er{e1a} have to be combinations of $e^{\pm
ik x}$.  The Wronskian $w$ for Case 3 is given by
\[
\begin{aligned}
\lb{dw1} w(k )=\{f_-(\cdot,k), f_+(\cdot,k )\}=ik f_+(0,k )+f_+'(0,k),
\end{aligned}
\]
where  $\{f, g\}=fg'-f'g$.  The functions $f_+(0,k), f_+'(0,k),
w(k)$ are entire, real on the imaginary line and  satisfy as $|k|\to
\iy, k\in \ol \C_+$:
\[
\lb{2.4}
\begin{aligned}
& f_+(0,k)=1+{O(1/k)},\qqq f_+'(0,k)=ik+O(1),\\
& w(k)=2ik+O(1),
\end{aligned}
\]
uniformly in $\arg \in [0,\pi]$. Let $F$ be one of the functions $f_+(0,\cdot), f_+'(0,\cdot)$ or
$w$.  The function $F$ has only finite simple  number of zeros in the upper half-plane $\C_+$  given by
\[
\lb{kn}
\begin{aligned}
&Case \ 1: \qqq\textstyle  i|E_1|^{1\/2},\dots ,  i|E_{m}|^{1\/2}\in i\R_+,
\\
&Case \ 2: \qqq \textstyle i|\wt E_1|^{1\/2},\dots , i|\wt
E_{N}|^{1\/2}\in i\R_+,
\\
&Case \ 3: \qqq \textstyle i|\cE_1|^{1\/2},\dots ,
i|\cE_{\cN}|^{1\/2}\in i\R_+,
\end{aligned}
\]
 and an infinite number of zeros, so-called
resonances, in  $\C_-$ and possibly one simple  zero at 0 (for $q\ne
0$). By definition, a zero of $F$ is called a resonance of the corresponding Schr\"odinger operator.  The
multiplicity of the resonance is the multiplicity of the
corresponding zero of $F$ and it can be any number, see \cite{K04}.
 Introduce the norm
$\|u\|^2=\int_0^1|u(x)|^2dx$ in $L^2(0,1)$.

%\newpage

%\newpage

%\section{Proof of the main results}
%\setcounter{equation}{0}
% Sect.4

\subsection{Proof Theorem \ref{T1}}
i-ii) {\bf Dirichlet  and Neumann  boundary conditions.}
Recall the identity from \cite{K04}:
\[
\lb{fs5} e^{-ik}f_+'(0,k)=i k\vt (1,k)-\vt'(1,k).
\]
From \er{fs5}, \er{aff1} and \er{bx} we have that $\wt E_1>\n_0$. Thus if $\n_0\ge 0$, then the operator
$\wt T$ has not eigenvalues.
We fix a real potential $q_0\in L^1(\R_+)$
 such that $\supp q_0\in [0,1]$ and $\n_0(q_0)=0$.
Define the potential $q_\ve, \ve\in \R$ by
$$
\supp q_\ve\ss[0,1], \qqq
q_\ve(x)=q_0(x)-\ve, \qqq x\in [0,1].
$$
We sometimes write $f_+(0,k,\ve), \m_{n}(\ve),..$ instead of $f_+(0,k), \m_{n},
..$ when several potentials $q_\ve$ are being dealt with.
The operators in \er{H01}-\er{H0110} and \er{defT} with the potential $q_\ve$ we denote by
$H_0(\ve), H_1(\ve),...  $ and  $T(\ve),...  $ and  the corresponding eigenvalues by $\m_n(\ve), \n_n(\ve)....$ and $E_n(\ve),....$.
Note that the eigenvalues $\m_n(\ve), \n_n(\ve),..$ satisfy
\[
\lb{ei2} \m_n(\ve)=\m_n(0)-\ve, \qq \n_n(\ve)=\n_n(0)-\ve,... \qq \t_n(\ve)=\t_n(0)-\ve, \qq
\vr_n(\ve)=\vr_n(0)-\ve.
\]
Firstly we consider the operator $\wt T$ with the {\bf  Neumann  boundary condition.}
Due to  \er{fs5}  the Jost function for  $\wt T(\ve)$ for $k=it, t\in \R$ satisfies
\[
\lb{ftt}
e^{-ik}f_+'(0,k,\ve)=-\vt'(1,k^2,\ve)+i k\vt (1,k^2,\ve)=i k\vt (1,\ve-t^2,0)-\vt'(1,\ve-t^2,0).
\]
 Define functions
$$
F(t,\ve)=-f(t^2,\ve)-tg(t^2,\ve),\qq f(t^2,\ve)=\vt'(1,\ve-t^2,q_0),\qq g(t^2,\ve)=\vt(1,\ve-t^2,q_0).
$$
 These functions are entire in $t,\ve $ and satisfy
\[
f(0,0)=0,\qqq  g(0,0)\ne 0,\qqq F(0,0)=0, \qq F_t(0,0)\ne 0.
\]
Then due to the Implicit Function Theorem there exists a
function $t(\ve)$, analytic in small disk $\{|\ve|<\t\}$ such that
$F(t(\ve),\ve)=0$ in the disk $\{|\ve|<\t\}$ (and here $\wt E_1(\ve)=-t^2(\ve)$ for $\ve>0$).

We have $t(\ve)=t_1\ve+O(\ve^2)$. Then from \er{ftt} we obtain for $t=t(\ve)$ as $\ve\to 0$:
\[
\lb{cxz}
\begin{aligned}
\vt'(1,\ve-t^2,0)=-t(\ve)\vt(1,\ve-t^2,0)=t_1\ve\vt(1,0,0)+O(\ve^2),
\\
\vt'(1,\ve -t^2,0)=\ve\dot \vt'(1,0,0)+O(\ve^2),
\end{aligned}
\]
which yields $t_1=-{\dot \vt'(1,0,0)\/\vt(1,0,0)}>0$, where $\dot u={\pa \/\pa \l}u$.

If $\ve>0$ is small enough, then $\n_0(\ve)=-\ve$ and
$\vr_1(\ve)=\vr_1(q_0)-\ve>0$, since there is the basic relation
\er{bx}. Thus we obtain $ \n_0(\ve)<\wt E_1(\ve)<0<\vr_1(\ve) $. If
$\ve$ is increasing then all eigenvalues $\n_0(\ve)<\wt E_1(\ve)<0<\vr_1(\ve)$ move monotonically to left and at $\ve_1=\t_1(0)$ we      have $\t_1({\ve_1})=0$.

Secondly, we consider the operator $T$ with the {\bf  Dirichlet boundary condition.}
The function $f_+(0,k)$ is expressed in terms of the fundamental
solutions $\vp ,\vt $  by
\[
\lb{fs4} e^{-ik}f_+(0,k)=\vp'(1,k^2)-ik\vp(1,k^2) \qq \forall \ k\in \C.
\]
Note that if the operator $H_{01}$ have eigenvalue $\t_1\ge 0$, then from \er{fs4}, \er{aff1} and \er{bx}
 we deduce that the operator $T$ has not any eigenvalue.
Let $\ve=\ve_1+z$, where $\ve_1=\t_1(0)$. Then due to \er{fs4} the Jost function $f_+(0,k,\ve)$ for the operator $T(\ve)$ with satisfies at $k=it, t\in \R$:
\[
\lb{fse}
e^{-ik}f_+(0,k,\ve)=\vp'(1,z-t^2,{\ve_1})+t\vp(1,z-t^2,{\ve_1}).
\]
We rewrite the rhs of the last identity in the form:
$$
F_o(t,\ve)=f_o(t^2,\ve)+tg_o(t^2,\ve), \qqq
f_o(t^2,z)=\vp'(1,z-t^2,{\ve_1}),\qqq g_o(t^2,z)=\vp(1,z-t^2,{\ve_1}).
$$
 These functions are entire in $t,z$ and satisfy
\[
f_o(0,0)=0, \qqq g_o(0,0)\ne 0,\qqq F_o(0,0)=0, \qq {\pa\/\pa t}F_o(0,0)\ne 0.
\]
Then the Implicit Function Theorem gives that there exists a
function $t_o(z)$, analytic in small disk $\{|z|<\d\}$ such that
$F(t_o(z),z)=0$ in the  disk $\{|z|<\d\}$. Here we have $E_1(\ve)=t_o^2(z), z=\ve -\ve_1$.
We have $t(z)=t_1z+O(z^2)$  as $z \to 0$. Then   from $f_+(0,k,\ve)=0$ and \er{fse} we obtain
\[
\lb{cxzx}
\begin{aligned}
\vp'(1,z-t^2,\ve_1)=-t(z)\vp(1,z-t^2,\ve_1)=
-t_1z\vp(1,0,\ve_1)+O(z^2),
\\
\vp'(1,z-t^2,\ve_1)=\dot\vp'(1,0,\ve_1)z+O(z^2),
\end{aligned}
\]
which yields $t_1=-{\dot\vp'(1,0,\ve_1)\/\vp(1,0,\ve_1)}>0$, since $\vp(1,\l,\ve_1)\to \iy$ as $\l\to -\iy$ and $\m_1>0$.

If $z>0$ is small enough, then $\t_1(\ve)=-z$ and
$\m_1(\ve)=\m_1(0)-\ve>0$, since there is the basic relation
\er{bx}. Thus we obtain
$$
\t_1(\ve)<E_1(\ve)<0<\m_1(\ve), \qq  {\rm and}\qq \wt E_1(\ve)<E_1(\ve).
$$
It is important that $ \wt E_1(\ve)<E_1(\ve)$ for any $\ve$ since $ \s_d(T)\cap \s_d(\wt T)=\es$.
If $\ve$ is increasing then all eigenvalues $\t_1(\ve)<E_1(\ve)<
\m_1(\ve)$ and $ \wt E_1(\ve)<E_1(\ve)$ move monotonically to left.
If $\ve$ is increasing more, then  we
have $\m_1(\ve_1)=0$ at $\ve_1=\m_1(0)$, but this eigenvalues does not "create" eigenvalues for the operators $T(\ve), \wt T(\ve)$.  If $\ve$ is increasing again  then  we have
$\n_1({\ve_2})=0$ at $\ve_2=\n_1(q_0)$. Using the above arguments for the
case $\n_0=0$ we obtain \er{D1} and \er{D2}, and so on.
Repeating these arguments we obtain \er{D1}-\er{N3}, since we can take any $q_0$ and $\ve$.

iii) We  consider the Schr\"odiger operator $\cT$ on {\bf the real line}.
Note that if the entire function  $ f_+'(0,k)<0$ for $k\in i\R_+$ and   has the zero $k=0$, then using
$w=ik f_+(0,\cdot)+f_+'(0,\cdot)$ from \er{dw1} we deduce that the operator $\cT$ has not any eigenvalue.
%From \er{dw1} and \er{N3} we deduce that  $\wt E_1 <\cE_1<E_1$.

  Recall that $q_0\in L^1(\R_+)$ is such that $\supp q_0\in [0,1]$ and $\n_0(q_0)=0$.
 We need   the identity from \cite{K05}:
\[
\lb{dw2}
e^{-ik}w(k )=2ik \D(\l)+\l\vp (1,\l)-\vt'(1,\l).
\]
 Due to \er{dw2} the Wronskian for  $q_\ve$ satisfies at $k=it$:
$$
\begin{aligned}
e^{-ik}w(k,\ve)= 2ik \D(\l,\ve)+\l\vp (1,\l,\ve)-\vt'(1,\l,\ve),
\\
e^{t}w(it,\ve)=-2t \D(\ve -t^2,0)-t^2\vp (1,\ve -t^2,0)-\vt'(1,\ve -t^2,0).
\end{aligned}
$$
We rewrite the rhs of the last identity  for  $k=it, t\in \R$ in the
form:
\[
\lb{F1}
F=-f-g, \
f(t,\ve)=\vt'(1,\ve -t^2,0),\  g(t,\ve)=2t \D(\ve -t^2,0)+t^2\vp (1,\ve -t^2,0).
\]
 These functions are entire in $t,\ve $ and satisfy
at $\ve=t=0$:
\[
\lb{Lf} f(0,0)=0,  \  f_t(0,0)=0,  \ g(0,0)=0, \qq g_t(0,0)=2\D(0,0)\ge 2,
\]
since $\D(\l,0)\ge 1$ for any $\l\le \n_0(q_0)$.
We show that exists exactly  one eigenvalue for small $\ve>0$.  From \er{dw2} we have the
simple fact:  $w(0,0)=0$ iff $\vt'(1,0,0)=0$.
Then the Implicit Function Theorem gives that there exists a
function $t_o(\ve)$, analytic in small disk $\{|\ve |<\d\}$ such that
$F(t_o(\ve),\ve)=0$ in the disk $\{|\ve|<\d\}$. Moreover, from \er{F1} we have $t_o(\ve)=-C\ve+O(\ve^2)$,
where $C={\dot \vt'(1,0,0)\/2\D(0,0)}<0$. Here we have $\cE_1(\ve)=-t_o^2(\ve)$.
Thus the eigenvalues $\n_0(\ve)$ at $\ve=0$ creates two eigenvalues $\n_0(\ve)<\wt E_1(\ve)<\cE_1(\ve)<0$ for small $\ve>0$. If $\ve$ is increasing then repeating the arguments for the case of the  operator $\wt T$ we have \er{N1x}, \er{N2x}. Here we need the fact: if $\cE$ is an eigenvalue of $\cT$, then  $\cE\notin \s_{d}(T)\cup  \s_{d}(\wt T)$.
\BBox

{\bf Proof of Corollary \ref{T2}}.  The Schr\"odinger equation
$
  -f''+q_e(x)f=k^2f,\ \  k\in \C \sm \{0\},
$ has the Jost solutions $\p_{\pm}(x,k)$  such that $\p_+(x,k
)=e^{ixk },\ \ x\ge 1$ and $\p_-(x,k)=e^{-ik x},\ \ x\leq -1$.
Note that the symmetry of the potential $q_e$ yields
$$
\p_+(x,k)=f_+(x,k)=\p_-(-x,k) \qqq \forall x\in [0,1].
$$
This implies that the Wronskian $w_e(k)$ for the potential $q_e$
satisfies
\[
\lb{tw} w_e(k)=\{\p_+(x,k),
\p_-(x,k)\}|_{x=0}=2f_+(0,k)f_+'(0,k),
\]
where $f_+(0,k)$ and $f_+'(0,k)$ are the Jost functions for the
Cases 1 and 2 respectively with the potential $q$. This identity and Theorem \ref{T1} yield
the proof of Corollary \ref{T2}.
\BBox

%\newpage

\

In order to discuss resonances  we define the class of all Jost
functions from \cite{K04}:

 {\it  {\bf Definition J.}  By   $\cJ$ we mean the class of all entire
functions $f$  having the form
\[
\begin{aligned}
f(k )=1+{\hat F(k)-\hat F(0)\/ 2ik}, \qqq \ k\in \C,
\end{aligned}
\]
where  $\hat F(k)=\int_0^1 F(x)e^{2ixk }dx$ is the Fourier
transformation  of $F\in\cP$ and a set of zeros $K=\{k_n, n\in \N\}$
(counted with multiplicity)  of $f$ satisfy:

1) The set $K=\ol K$ and has not zeros from $\R\sm \{0\}$ and has
possibly one simple zero at 0.

2) Let  $\# E$ be  the number of zeros of $K$ on the set $E\ss \C$.
The set $K$  has finite number elements $k_1,..,k_m$ from $\C_+$,
which are simple, belong to $i\R_+$ and if they are labeled by
$|k_1|>|k_2|>...>|k_m|>0$, then  intervals on $i\R_-$ defined
  $I_j=(-k_j,-k_{j+1}), j\in \N_{m-1}$ and $I_m=(-k_m, 0]$ satisfy}
\[
\begin{aligned}
\lb{c} -k_j\notin K,\ \  j\in \N_m:=\{1,2,...,m\},
% (-1)^jf(-k_j) >0,\ \qqq j\in \N_m:=\{1,2,...,m\},
\qq \ca \# I_j\ge 1 \qq {\rm is\ odd},\  j\in \N_{m-1}
\\
  \# I_m    \ge 0 \qq {\rm is\ even}              \ac.
\end{aligned}
\]

%\er{c} is equivalent  to  the following   condition: $f(-k_j)\neq
%0,\ \ j\in \N_m$ and

%Let  $\#(f,E)$ be  the number of zeros of $f$ (counted with
%multiplicity) in the set $E\ss \C$. }

%Define the intervals $I_j=(-k_{j},-k_{j+1}), j\in \N_{m-1},$ and
%$I_m=(-k_{m},0]$ on by $i\R_-$. Let  $\#(E)$ be  the number of zeros
%of $\p$ (counted with multiplicity) in the set $E\ss \C$.

%The zeros of function $\p$ on $I_j$ satisfies:
%\[ \#(I_j)=\gn_j, \qq \ca \gn_j\ge 1 \ {\rm is \ odd \ if} \ j\in \N_{m-1}\\
%                     \gn_m\ge 0 \ {\rm is \ even}\ac \]

We recall the basic result about inverse problem (including the
characterization) for compactly supported potentials from
\cite{K04}.

\begin{theorem}
\lb{TA}
 The mapping $\p: \cP\to \cJ$ given by $q\to  \p$ is one-to-one and onto,
 where $\p(k)=f_+(0,k)$.
\end{theorem}

Thus if $q\in \cP$, then $\p\in \cJ$. Conversely,  for each $f\in
\cJ$ exists a unique potential $q\in \cP$ such that the Jost
function $\p=f$. In fact using this result we reformulate the
problem for the differential operator as the problem of the entire
function theory. We need Theorem 3 from \cite{K04x}:

\begin{theorem}
\lb{TB} Let $\p^o\in \cJ$ have zeros $\{k_n^o\}_{1}^\iy$ and let the
sequence complex points $\{k_n\}_{1}^\iy$ satisfy 1) and 2) in
Condition J and  $\sum_{n\ge 1} n^3|k_n-k_n^o|^2 <\iy$. Then $\{
k_n\}_{1}^\iy$ is a sequence of zeroes of unique $\p\in \cJ$, which
is the Jost function for unique $q\in \cP$.
\end{theorem}

{\bf Proof of Theorem \ref{T1x}}. i) Let a resonance   $k_o=-ir\in
I_j$. Then from \er{fs4} we obtain $\vp'(1,k_o^2)=r\vp(1,k_o^2)$,
which yields $\sign\vp'(1,k_o^2)=\sign \vp(1,k_o^2)$. From \er{D2}
we deduce that $\m_j<r_o^2<\t_{j+1}$, since $\vp'(1,\l)\to +\iy$ and
$\vp(1,\l)\to +\iy$ as $\l\to -\iy$.

%which is the Jost function for unique $q\in \cP$.

ii) We take any potential $q^o\in \cP$ such that the Jost function
$\p^o\in \cP$ has zeros $\{k_n^o\}_{1}^\iy$ only from $\C_-$. We
take the half-disc $\cD_r:=\{k\in \C_-: |k|<r \}$ for a radius
$r>3+|p_1|$ large enough such that the number $N_r$ of zeros
$\{k_n^o\}$ in $\cD_r$ satisfy $N_r>m+(s_1+..+s_m)$, since we have
the following result from \cite{Z87}: $N_r={r\/2\pi}(1+o(1)$ ar
$r\to \iy$. We construct the new sequence of zeros $\{k_n\}$:

$\bu $ if $k_n^o\notin \cD_r$, then $k_n=k_n^o$.

$\bu $ Consider $N_r$ zeros $\{k_n^o\}$ in  $\cD_r$. We have three
cases:

a) We remove $m$ zeros from $\dD_-(r)$  on the points $k_j=p_j\in
i\R_+$ for each $j\in \N_m$.

b) We remove $s_j$ zeros from $\dD_-(r)$ on each interval
$I_j=(-k_{j},-k_{j+1}), j\in \N_{m-1}$ and $I_m=(-k_{m},0]$ at
$j=m$.

c) We remove all remaining zeros $n_r=N_r-m-(s_1+..+s_m)>0$ from
$\dD_-(r)$ on point $-p_1-i$ with the multiplicity $n_r$. Recall
that $r$ is large enough and $r>3+|p_1|$.

Thus the sequence $\{k_n\}$ satisfy 1) and 2) in Condition J. Then
due to Theorem \ref{TB} $(k_n)_1^\iy$ is a sequence of zeroes of
unique $f\in \cJ$,  which is the Jost function for unique $q\in
\cP$. \BBox

%\newpage

\no {\bf Proof of Corollary \ref{T3}}. i) We need an inequality from
\cite{HLT98}, \cite{W96}: if $q\in L^1(\R)$, then
 \[
\lb{W} \sum_{\cE_n<0}|\cE_n|^{1\/2}\le {1\/2}\int_\R q_-(x)dx.
\]
From Theorem \ref{T1} and inequalities \er{Be}, \er{CC} and \er{W}
we obtain \er{nu1}, \er{nu3} and \er{nu2} for
$\gS=\sum_{\m_n<0}|\m_n|^{1\/2}$. Applying \er{W} to the operator
$\cT_e$ and estimates from \er{bx},  Corollary \ref{T2}
($|\l_n^\pm|\le |\cE_n^e|$ for negative eigenvalues) we obtain
$$
 \sum_{n>0,
\l_n^\pm<0}|\l_n^\pm|^{1\/2}<\sum_{\cE_n^e<0}|\cE_n^e|^{1\/2}
 \le {1\/2}\int_\R q_{e-}(x)dx= \int_0^1 q_-(x)dx.
$$
ii) Substituting the estimate $\|q\|\le 2\g\max \{1, \g^{1\/3}\}$
from \er{esg1} into \er{nu1}, \er{nu2}  we obtain \er{nu2x}. \BBox

%\newpage

We now illustrate  our results

{\bf Example}\lb{T6} {\it i) For any number $c<0$ there exists a
potential $q\in \cP$ such that the operator $T_q$ has not
eigenvalues and $\int_0^1q(x)dx=c$.
%\[ \lb{aa1} \s(T_q)=\s_{ac}(T_q),\qqq \int_0^1q(x)dx=q_0. \]

ii) Let $c>0$ and let $I_j=(t_j,m_j)\ss \R_-, j\in \N_m$ be   finite
sequence of intervals such that $ t_1<m_1<t_2<m_2<.... <t_m<m_m<0$.
Then there exists a potential $q\in \cP$ such that
$\int_0^1q(x)dx=c$ and the operator $T_q$ has an eigenvalue $E_j\in
I_j$ for any $j\in \N_m$. There are no other eigenvalues. }

{\bf Proof.} i) Let $c<0$.  We define two strongly increasing
sequences by
$$
%\begin{aligned}
m_j=\ca (\pi j)^2, \qqq  & j\in \N_n \\ (\pi j)^2+c, \qq & j>n\ac,
\qqq \textstyle t_j=\ca \pi^2 (j-{1\/2})^2, & j\in \N_n \\
\pi^2 (j-{1\/2})^2+c, & j>n \ac,
%\end{aligned}
$$
for some $n>|c|$. They satisfy  $t_j<m_j<t_{j+1}<... j\in \N$. Due
to result of \cite{M86} there exists a potential $q\in L^2(0,1)$
with $\int_0^1q(x)dx=c$ and $\m_j(q)=m_j$ and $\t_j(q)=t_j$ for all
$j\ge 1$. Then due to  \er{D2} the operator $T=-{d^2\/dx^2}+q$ on
$\R_+$ has not eigenvalues.

ii)  Let $c>0$. Let $t_j=\pi^2 (j-{1\/2})^2+c, j>m$ and $m_j=(\pi
j)^2+c, j>m$. Thus the numbers $t_j,m_j, j\in \N$ satisfy
$t_j<m_j<t_{j+1}$ for all $j\in \N$. Due to result of \cite{M86}
there exists a potential $q\in L^2(0,1)$ with $\int_0^1q(x)dx=c$ and
$\m_j(q)=m_j$ and $\t_j(q)=t_j$ for all $j\ge 1$. Then due to
\er{D2} the operator $T=-{d^2\/dx^2}+q$ on $\R_+$ has an eigenvalue
$E_j\in I_j$ for any $j\in \N_m$ and there are no other eigenvalues.
 \BBox

\

\footnotesize \no {\bf Acknowledgments.} \footnotesize This work was
supported by the RSF grant  No. 18-11-00032.

\end{document}